\documentclass[preprint,3p,12pt]{elsarticle}

\usepackage{lineno,hyperref}
\modulolinenumbers[5]

\usepackage{hyperref}

\usepackage{epstopdf}

\usepackage{epigraph}

\usepackage{mathrsfs}
\usepackage{amssymb}
\usepackage{amsfonts}
\usepackage{latexsym}
\usepackage{amsmath}
\usepackage{xcolor}
\usepackage{wrapfig}
\usepackage{floatflt}

\usepackage{graphicx}
\usepackage{subcaption}
\def\lb{\label}

\newcommand{\er}[1]{\textrm{(\ref{#1})}}

\newtheorem{theorem}{\bf Theorem}[section]

\def\a{\alpha}

\def\d{\delta}

\def\m{\mu}

\def\o{\omega}

\def\Z{{\mathbb Z}}       \def\C{{\mathbb C}}

\def\lt{\biggl}                  \def\rt{\biggr}

\let\ge\geqslant                 \let\le\leqslant

\def\iy{\infty}
               
                 \def\ts{\times}
\def\pa{\partial}

\def\el2{\ell^{\,2}}             \def\1{1\!\!1}

\newtheorem{corollary}[theorem]{\bf Corollary}

\let\ge\geqslant
\let\le\leqslant

\newcommand{\ca}{\begin{cases}}
\newcommand{\ac}{\end{cases}}
\newcommand{\ma}{\begin{pmatrix}}
\newcommand{\am}{\end{pmatrix}}
\def\eq{\begin{equation}}
\def\qe{\end{equation}}
\renewcommand{\[}{\begin{equation}}
\renewcommand{\]}{\end{equation}}

\begin{document}

\begin{frontmatter}

\title{A note on exotic integrals}
\date{\today}

\author
{Anton A. Kutsenko}

\address{Jacobs University (International University Bremen), 28759 Bremen, Germany; email: akucenko@gmail.com}

\epigraph{In Memory of Robert S. Strichartz (1943-2021)}

\begin{abstract}
 We consider Bernoulli measures $\mu_p$ on the interval $[0,1]$. For the standard Lebesgue measure the digits $0$ and $1$ in the binary representation of real numbers appear with an equal probability $1/2$. For the Bernoulli measures, the digits $0$ and $1$ appear with probabilities $p$ and $1-p$, respectively. We provide explicit expressions for various $\mu_p$-integrals. In particular, integrals of polynomials are expressed in terms of the determinants of special Hessenberg matrices, which, in turn, are constructed from the Pascal matrices of binomial coefficients. This allows us to find closed-form expressions for the Fourier coefficients of $\mu_p$ in the Legendre polynomial basis. At the same time, the trigonometric Fourier coefficients are values of some special entire function, which admits an explicit infinite product expansion and satisfies interesting properties, including connections with the Stirling numbers and the polylogarithm.
\end{abstract}

\begin{keyword}
Bernoulli measure on intervals
\end{keyword}


\end{frontmatter}


{\section{Introduction}\lb{sec1}}

This short note is motivated by a nice paper \cite{S}, where the Bernoulli measure $\m_p$ is defined as an ``exotic" modification of the Lebesgue measure. Roughly speaking, the Lebesgue measure means that each digit in the binary representation of a real number appears with the probability $\frac12$. For Bernoulli measures, digit $0$ appears with the probability $p$, and digit $1$ appears with the probability $1-p$. For most of $p$, the resulting Bernoulli measure has a complex singular fractal structure. In \cite{S}, recurrent formulas for $\int_0^1 x^n d\m_p(x)$ are presented. In the current paper, we express $\int_0^1 {Polynomial}(x) d\m_p(x)$  through determinants of special Hessenberg matrices, which, in turn, are perturbations of Pascal matrices. This representation can be also convenient for the expression of $\m_p$-integrals of transcendent functions.

The Bernoulli measure $\m_p$, $p\in(0,1)$ is a probabilistic Borel measure on $[0,1]$ that satisfies
$$
 \m_p([0,\frac12])=p,\ \ \m_p([\frac12,1])=1-p,
$$
$$
 \m_p([0,\frac14])=p^2,\ \ \m_p([\frac14,\frac12])=p(1-p),\ \ \m_p([\frac12,\frac34])=(1-p)p,\ \ \m_p([\frac34,1])=(1-p)^2
$$
and so on. The general formula for this measure on dyadic intervals is
\[\lb{001}
 \m_p(\{0.x_1x_2...x_ny_1y_2....:\ y_i\in\{0,1\}\})=\prod_{i=1}^np^{1-x_i}(1-p)^{x_i},\ \ \ x_i\in\{0,1\},
\]
where $0.x_1x_2...x_ny_1y_2....$ denotes the $2$-base digital expansion of a real number from $[0,1]$. Let $f$ be some continuous function. It is not difficult to check that
\[\lb{002}
 \int_0^1 f(x)d\m_p(x)=p\int_0^1f\lt(\frac x2\rt)d\m_p(x)+(1-p)\int_0^1f\lt(\frac{1+x}2\rt)d\m_p(x),
\]
see \cite{S}. If $p=\frac12$ then $\m_p$ is the standard Lebesgue measure. For the limit cases $p=1$ and $p=0$, $\m_p$ is the Dirac $\d$-measure located at $x=0$ and $x=1$ respectively.

 We denote $|{\bf A}|$  the determinant of square matrix ${\bf A}$. We also use the notation $\binom{n}{m}$ for the binomial coefficients, including the generalized ones, see identity \er{010} in Corollary \ref{C1}. Let us formulate the main result.

\begin{theorem}\lb{T1}
i) Let $R(x)=\sum_{n=0}^N r_n x^n$ be some polynomial with $r_n\in\C$. Then
\[\lb{003}
 \int_0^1 R(x)d\m_p(x)=\frac{(1-p)^N}{\prod_{n=1}^N(2^n-1)}\begin{vmatrix} 
     1 & \frac{1-2}{1-p} & 0 & ... &  0 \\
     1 & \binom{2}{1} & \frac{1-2^2}{1-p} & ... &  0 \\
     ... & ... & ... & ... & ...  \\
     1 & \binom{N}{1} & \binom{N}{2} & ... & \frac{1-2^N}{1-p}  \\
     r_0 & r_1 & r_2 & ... & r_N 
 \end{vmatrix}.
\]
In particular,
\[\lb{004}
 \int_0^1 x^Nd\m_p(x)=\frac{(1-p)^N}{\prod_{n=1}^N(2^n-1)}\begin{vmatrix} 
     1 & \frac{1-2}{1-p} & 0 & ... &  0 \\
     1 & \binom{2}{1} & \frac{1-2^2}{1-p} & ... &  0 \\
      1 & \binom{3}{1} & \binom{3}{2} & ... &  0 \\
     ... & ... & ... & ... & ...  \\
     1 & \binom{N}{1} & \binom{N}{2} & ... & \binom{N}{N-1}  
 \end{vmatrix}.
\]
ii) Let us denote $I_p(\o)=\int_0^1e^{\o x}d\m_p(x)$ for $\o\in\C$. Then $I_p(\o)$ is an entire function, which satisfies the functional equation
\[\lb{005}
 I_p(\o)=(p+(1-p)e^{\frac{\o}2})I_p\lt(\frac{\o}2\rt),
\]
that, in turn, leads to
\[\lb{006}
 I_p(\o)=\prod_{n=1}^{+\iy}(p+(1-p)e^{\frac{\o}{2^n}}).
\]
Product \er{006} converges absolutely and uniformly on any compact subset of $\C$. Taylor series for $I_p(\o)$ can be written explicitly
\[\lb{007}
 I_p(\o)=1+(1-p)\o+\sum_{N=2}^{+\iy}\frac{(1-p)^N\o^N}{N!\prod_{n=1}^N(2^n-1)}\begin{vmatrix} 
     1 & \frac{1-2}{1-p} & 0 & ... &  0 \\
     1 & \binom{2}{1} & \frac{1-2^2}{1-p} & ... &  0 \\
      1 & \binom{3}{1} & \binom{3}{2} & ... &  0 \\
     ... & ... & ... & ... & ...  \\
     1 & \binom{N}{1} & \binom{N}{2} & ... & \binom{N}{N-1}  
 \end{vmatrix}.
\]
Function $I_p(\o)$ admits analytic continuation to all $(p,\o)\in\C\ts\C$. Moreover, it satisfies the symmetry and reciprocal properties
\[\lb{008}
I_p(-\o)=e^{-\o}I_{1-p}(\o),\ \ \ \ \ \prod_{m=1}^{+\iy}I_{e^{\frac{p}{2^m}}}(\o)=\prod_{m=1}^{+\iy}I_{e^{\frac{\o}{2^m}}}(p).
\]
If $p\ne1/2$ then there is an additional property
\[\lb{008a}
 I_{p}\lt(\frac{\o}2\rt)I_{\frac{p}{2p-1}}\lt(\frac{\o}2\rt)=I_{\frac{p^2}{2p-1}}(\o).
\]
There is also a notable series for the logarithm 
\[\lb{new1}
 \ln I_p(\o)=\sum_{n=1}^{+\iy}Q_n(p)\frac{\o^n}{(2^{n}-1)n!},
\] 
where the polynomials $Q_n(p)$ can be defined as
\[\lb{new2}
 Q_1(p)=1-p,\ \ \ Q_{n+1}(p)=p(p-1)Q_n'(p),\ \ \ n\ge1.
\]
All zeros of $Q_n(p)$ are simple and belong to $[0,1]$. The series \er{new1} converges for $|\o|<2|\ln\frac{p-1}{p}|$. Further analysis allows us to express $Q_n(p)$ through the polylogarithm of non-positive integer orders
\[\lb{new3}
 Q_n(p)=-{\rm Li}_{1-n}\lt(\frac{p-1}{p}\rt)=-\sum_{k=0}^{n-1}k!S(n,k+1)(p-1)^{k+1},
\]
where the Stirling numbers of the second kind are given by
\[\lb{new4}
 S(n,k)=\frac1{k!}\sum_{i=0}^k(-1)^i\binom{k}{i}(k-i)^n.
\]
\end{theorem}

{\bf Remark 1.} Combining \er{005} and \er{008a}, we obtain
\[\lb{009a}
 I_p(\o)=\frac{p+(1-p)e^{\frac{\o}2}}{I_{\frac{p}{2p-1}}(\frac{\o}2)}I_{\frac{p^2}{2p-1}}(\o),
\]
if the denominator is non-zero. Formula \er{009a} is an analog of \er{005} for the first argument, since it allows to reduce it twice:
$$
 \frac{p^2}{2p-1}=\frac{p}{2}+O(1)
$$
for large $p$. Moreover, we have
$$
 \frac{p}{2p-1}=\frac12+O\lt(\frac1p\rt)
$$
for large $p$. Thus, \er{009a} can be useful for numerical computations when $p$ is large. Applying \er{009a} to its second term in RHS many times, one may obtain an analog of \er{006}, containing a nonlinear iteration of arguments. As an interesting exercise, using \er{009a}, one may obtain a correct order of growth of $I_p(\o)$ as $e^{\frac{\ln^2 p}{\ln 4}+O(\ln p)}$ for $p\to+\iy$.

{\bf Remark 2.} For $k\in\Z$, identity \er{005} implies
$$
 I_p(2\pi i \cdot k)=I_p(2\pi i \cdot 2k)=I_p(2\pi i \cdot 2^2k)=I_p(2\pi i \cdot2^3 k)=....
$$
Thus, the even trigonometric Fourier coefficients of $\m_p$ coincides with some odd ones $I_p(2\pi(2n+1))$, $n\in\Z$. The Fourier coefficients of $\m_p$ in the shifted Legendre polynomial system can be computed explicitly
$$
 \int_0^1P_N(x)d\m_p(x)=\frac{(p-1)^N}{\sqrt{2N+1}\prod_{n=1}^N(2^n-1)}\begin{vmatrix} 
     1 & \frac{1-2}{1-p} & 0 & ... &  0 \\
     1 & \binom{2}{1} & \frac{1-2^2}{1-p} & ... &  0 \\
     ... & ... & ... & ... & ...  \\
     1 & \binom{N}{1} & \binom{N}{2} & ... & \frac{1-2^N}{1-p}  \\
     1 & -\binom{N}{1}\binom{N+1}{1} & \binom{N}{2}\binom{N+2}{2} & ... & (-1)^N\binom{N}{N}\binom{2N}{N}
 \end{vmatrix},
$$  
where $P_N(x)=\frac{(-1)^N}{\sqrt{2N+1}}\sum_{n=0}^N\binom{N}{n}\binom{N+n}{n}(-x)^n$.

Finally, using Taylor expansions of various transcendent functions, we may apply \er{003} to obtain the following identities for \er{006} and other integrals.

\begin{corollary}\lb{C1}
The following identities are fulfilled 
\[\lb{009}
 I_p(\o)=\int_{0}^1e^{\o x}d\m_p(x)=\lim_{N\to\iy}\frac{(1-p)^N}{\prod_{n=1}^N(2^n-1)}\begin{vmatrix} 
     1 & \frac{1-2}{1-p} & 0 & ... &  0 \\
     1 & \binom{2}{1} & \frac{1-2^2}{1-p} & ... &  0 \\
     ... & ... & ... & ... & ...  \\
     1 & \binom{N}{1} & \binom{N}{2} & ... & \frac{1-2^N}{1-p}  \\
     1 & \o & \frac{\o^2}{2!} & ... & \frac{\o^N}{N!} 
 \end{vmatrix},\ \ \ \o\in\C,
\]
\[\lb{010}
 \int_{0}^1\ln(1-\o x)d\m_p(x)=\lim_{N\to\iy}\frac{-(1-p)^N}{\prod_{n=1}^N(2^n-1)}\begin{vmatrix} 
     1 & \frac{1-2}{1-p} & 0 & ... &  0 \\
     1 & \binom{2}{1} & \frac{1-2^2}{1-p} & ... &  0 \\
     ... & ... & ... & ... & ...  \\
     1 & \binom{N}{1} & \binom{N}{2} & ... & \frac{1-2^N}{1-p}  \\
     0 & \o & \frac{\o^2}2 & ... & \frac{\o^N}N 
 \end{vmatrix},\ \ \ |\o|<1,
\]
\[\lb{011}
 \int_{0}^1(1+\o x)^{\a}d\m_p(x)=\lim_{N\to\iy}\frac{(1-p)^N}{\prod_{n=1}^N(2^n-1)}\begin{vmatrix} 
     1 & \frac{1-2}{1-p} & 0 & ... &  0 \\
     1 & \binom{2}{1} & \frac{1-2^2}{1-p} & ... &  0 \\
     ... & ... & ... & ... & ...  \\
     1 & \binom{N}{1} & \binom{N}{2} & ... & \frac{1-2^N}{1-p}  \\
     1 & \binom{\a}{1}\o & \binom{\a}{2}\o^2 & ... & \binom{\a}{N}\o^N 
 \end{vmatrix},\ \ \ |\o|<1.
\]
Moreover, RHS of \er{009}-\er{011} is exact starting segment in Taylor expansion of analytic functions of $\o$ represented by LHS of \er{009}-\er{011}, respectively.
\end{corollary}

{\section{Proof of Theorem \ref{T1}}\lb{sec2}}

Let us denote 
\[\lb{100}
 J_n:=\int_0^1x^nd\m_p(x),\ \ n\ge0.
\] 
We fix some $N\ge1$. Using \er{002}, we obtain
\[\lb{101}
 p2^{-N}J_N+(1-p)2^{-N}\sum_{n=0}^N\binom{N}{n}J_n=J_N,
\]
that can be written as
\[\lb{102}
 (1-p)2^{-N}\sum_{n=0}^{N-1}\binom{N}{n}J_n+(2^{-N}-1)J_N=0
\]
or in the matrix form
\[\lb{103}
 \ma 1 & 0 & 0 & ... & 0 \\
     \frac{(1-p)\binom{1}{0}}2 & \frac12-1 & 0 & ... & 0 \\
     \frac{(1-p)\binom{2}{0}}{2^2} & \frac{(1-p)\binom{2}{1}}{2^2} & \frac1{2^2}-1 & ... & 0 \\
     ... & ... & ... & ... & ... \\
     \frac{(1-p)\binom{N}{0}}{2^N} & \frac{(1-p)\binom{N}{1}}{2^N} & \frac{(1-p)\binom{N}{2}}{2^N} & ... & \frac1{2^N}-1
 \am\ma J_0 \\ J_1 \\ J_2 \\ ... \\ J_N \am=\ma 1 \\ 0 \\ 0 \\ ... \\ 0 \am.
\]
Using the Cramer's rule, we write
\begin{multline}\lb{104}
 \int_0^1 R(x)d\m_p(x)=\sum_{n=0}^Nr_nJ_n=\frac{-1}{\prod_{n=1}^N(2^{-n}-1)}\begin{vmatrix} 
  1 & 0 & 0 & ... & 0 & 1 \\
     \frac{(1-p)\binom{1}{0}}2 & \frac12-1 & 0 & ... & 0 & 0 \\
     \frac{(1-p)\binom{2}{0}}{2^2} & \frac{(1-p)\binom{2}{1}}{2^2} & \frac1{2^2}-1 & ... & 0 & 0 \\
     ... & ... & ... & ... & ... & ... \\
     \frac{(1-p)\binom{N}{0}}{2^N} & \frac{(1-p)\binom{N}{1}}{2^N} & \frac{(1-p)\binom{N}{2}}{2^N} & ... & \frac1{2^N}-1 & 0 \\
     r_0 & r_1 & r_2 & ... & r_N & 0
 \end{vmatrix}=\\
 \frac{1}{\prod_{n=1}^N(1-2^{-n})}\begin{vmatrix} 
     \frac{(1-p)\binom{1}{0}}2 & \frac12-1 & 0 & ... &  0 \\
     \frac{(1-p)\binom{2}{0}}{2^2} & \frac{(1-p)\binom{2}{1}}{2^2} & \frac1{2^2}-1 & ... &  0 \\
     ... & ... & ... & ... & ...  \\
     \frac{(1-p)\binom{N}{0}}{2^N} & \frac{(1-p)\binom{N}{1}}{2^N} & \frac{(1-p)\binom{N}{2}}{2^N} & ... & \frac1{2^N}-1  \\
     r_0 & r_1 & r_2 & ... & r_N 
 \end{vmatrix},
\end{multline}
that, after elementary transformations, leads to \er{003}. Note that, in some sense, the last identity in \er{104} is more convenient than \er{003}, since, e.g., there is no singularity at $p=1$. We place \er{003} in the statement of the theorem because it looks more compact. Using \er{002}, we get
$$
 I_p(\o)=\int_0^1 e^{\o x}d\m_p(x)=p\int_0^1e^{\frac{\o x}2}d\m_p(x)+(1-p)\int_0^1e^{\frac{\o (x+1)}2}d\m_p(x)=(p+(1-p)e^{\frac{\o}2})I_p\lt(\frac{\o}2\rt),
$$
that gives \er{005}. Repeating \er{005} and using $I_p(0)=1$ we arrive at \er{006}. The terms of \er{006} have the form
$$
 p+(1-p)e^{\frac{\o}{2^n}}=1+(p-1)(e^{\frac{\o}{2^n}}-1),\ \ \ |(p-1)(e^{\frac{\o}{2^n}}-1)|\le 2^{-n}C(p,\o),
$$
where we can take 
$$
 C(p,\o)=|p-1|(e^{|\o|}-1),
$$ 
which is bounded on compact subsets of $\C\ts\C$ and does not depend on $n$. Hence, the Weierstrass criterion is satisfied, and $I_p(\o)$ is entire in both variables $p,\o\in\C$. Using \er{006}, we obtain the following properties
$$
 I_p(-\o)=\prod_{n=1}^{+\iy}(p+(1-p)e^{\frac{-\o}{2^n}})=\prod_{n=1}^{+\iy}(pe^{\frac{\o}{2^n}}+(1-p))e^{\frac{-\o}{2^n}}=e^{-\o}I_{1-p}(\o)
$$
and
$$
 \prod_{m=1}^{+\iy}I_{e^{\frac{p}{2^m}}}(\o)=\prod_{m=1}^{+\iy}\prod_{n=1}^{+\iy}(e^{\frac{p}{2^m}}+(1-e^{\frac{p}{2^m}})e^{\frac{\o}{2^n}})=\prod_{m=1}^{+\iy}\prod_{n=1}^{+\iy}(1+(e^{\frac{p}{2^m}}-1)(e^{\frac{\o}{2^n}}-1))=\prod_{m=1}^{+\iy}I_{e^{\frac{\o}{2^m}}}(p)
$$
which leads to \er{008}. Using \er{006} again, we obtain also
\begin{multline}\notag
 I_{p}\lt(\frac{\o}2\rt)I_{\frac{p}{2p-1}}\lt(\frac{\o}2\rt)=\prod_{n=1}^{+\iy}(p+(1-p)e^{\frac{\o}{2^{n+1}}})\prod_{n=1}^{+\iy}(\frac{p}{2p-1}+(1-\frac{p}{2p-1})e^{\frac{\o}{2^{n+1}}})=\\
 \prod_{n=1}^{+\iy}(p+(1-p)e^{\frac{\o}{2^{n+1}}})(\frac{p}{2p-1}+(1-\frac{p}{2p-1})e^{\frac{\o}{2^{n+1}}})=\prod_{n=1}^{+\iy}(p+(1-p)e^{\frac{\o}{2^{n+1}}})(\frac{p}{2p-1}+\frac{p-1}{2p-1}e^{\frac{\o}{2^{n+1}}})=\\
 \prod_{n=1}^{+\iy}(\frac{p^2}{2p-1}+(1-\frac{p^2}{2p-1})e^{\frac{\o}{2^n}})=I_{\frac{p^2}{2p-1}}(\o),
\end{multline}
which gives \er{008a}. Let us denote
\[\lb{105}
 J(p,\o)=\ln(p+(1-p)e^{\o})=\sum_{n=1}^{+\iy}Q_n(p)\frac{\o^n}{n!},
\]
where the radius of convergence of the series coincide with the norm of the smallest zero of $p+(1-p)e^{\o}$, i.e. with $|\ln\frac{p}{p-1}|$. Using \er{105}, we obtain 
\[\lb{106}
 p(p-1)\frac{\pa J}{\pa p}(p,\o)=\frac{p(p-1)(1-e^{\o})}{p+(1-p)e^{\o}}=p-\frac{p}{p+(1-p)e^{\o}}
\]
and
\[\lb{107}
 \frac{\pa J}{\pa \o}(p,\o)=1-\frac{p}{p+(1-p)e^{\o}},
\]
which gives
\[\lb{108}
 p(p-1)\frac{\pa J}{\pa p}(p,\o)+1-p=\frac{\pa J}{\pa \o}(p,\o),
\]
and, hence, we arrive at \er{new2}.
According to the Rolle's Theorem and the principle of mathematical induction, all zeros of $Q_n(p)$ are simple and belong to $[0,1]$, see \er{new2}. Using \er{005} and \er{106}, we obtain
\[\lb{110}
 \ln I_p(\o)=\sum_{n=1}^{+\iy}J\lt(p,\frac{\o}{2^n}\rt)=\sum_{n=1}^{+\iy}\sum_{m=1}^{+\iy}Q_m(p)\frac{\o^m}{2^{nm}m!}=\sum_{m=1}^{+\iy}Q_m(p)\frac{\o^m}{(2^{m}-1)m!},
\]
which leads to \er{new1}. Taking some fixed $\o$ and sufficiently small $1-p$, we expand \er{105} as
\begin{multline}\lb{111}
 \ln(p+(1-p)e^{\o})=\ln p+\ln\lt(1-\frac{(p-1)e^{\o}}{p}\rt)=\ln p-\sum_{m=1}^{+\iy}\frac{(p-1)^me^{m\o}}{mp^m}=\ln p+\\
 -\sum_{m=1}^{+\iy}\frac{(p-1)^m}{mp^m}\sum_{n=0}^{+\iy}\frac{(m\o)^{n}}{n!}=\ln p-\sum_{m=1}^{+\iy}\frac1m\lt(\frac{p-1}{p}\rt)^m-\sum_{n=1}^{+\iy}\frac{\o^n}{n!}\sum_{m=1}^{+\iy}m^{n-1}\lt(\frac{p-1}{p}\rt)^m=\\
 \ln p+\ln\lt(1-\frac{p-1}{p}\rt)-\sum_{n=1}^{+\iy}\frac{\o^n}{n!}{\rm Li}_{1-n}\lt(\frac{p-1}{p}\rt)=-\sum_{n=1}^{+\iy}\frac{\o^n}{n!}{\rm Li}_{1-n}\lt(\frac{p-1}{p}\rt),
\end{multline}
where the polylogarithm is defined by
\[\lb{112}
 {\rm Li}_s(z)=\sum_{m=1}^{+\iy}\frac{z^m}{m^s}.
\]
Comparing \er{111} with \er{105} we deduce the first identity in \er{new3}. The second one follows from the known representation of the polylogarithm of non-positive integer orders $s$ through the Stirling numbers, see, e.g., \cite{L}
\[\lb{113}
 {\rm Li}_{-n}(z)=\sum_{k=0}^nk!S(n+1,k+1)\lt(\frac{z}{1-z}\rt)^{k+1},\ \ \ n=0,1,2,3,...,
\] 
where the Stirling numbers $S(n,k)$ are given by \er{new4}.

\end{document}